\documentstyle[amscd,amssymb,verbatim,12pt]{amsart}
\newcommand{\bx}{{\bf x}}
\newcommand{\bv}{{\bf v}}
\newcommand{\bw}{{\bf w}}
\newcommand{\bu}{{\bf u}}
\newcommand{\by}{{\bf y}}
\newcommand{\be}{{\bf e}}
\newcommand{\bd}{{\bf d}}
\newcommand{\bc}{{\bf c}}
\newcommand{\bla}{{\bf \lambda}}
\newcommand{\Bases}{\operatorname{Bases}}

\renewcommand{\mod}{\operatorname{mod}}

\newcommand{\OO}{{\cal O}}

\newcommand{\SL}{\operatorname{SL}}
\newcommand{\G}{{\Bbb G}}

\newcommand{\lan}{\langle}
\newcommand{\ran}{\rangle}

\newcommand{\gH}{{\frak H}}

\newcommand{\Th}{\Theta}

\newcommand{\Pic}{\operatorname{Pic}}

\renewcommand{\ker}{\operatorname{ker}}

\numberwithin{equation}{section}

\newtheorem{thm}{Theorem}[section]

\newtheorem{prop-def}[thm]{Proposition-Definition}
\newtheorem{lem}[thm]{Lemma}

\newenvironment{defi}{\vspace{3mm}\noindent
{\bf Definition.}}{\vspace{3mm}}

\newcommand{\Pf}{\noindent {\it Proof}}
\newcommand{\id}{\operatorname{id}}

\newcommand{\ov}{\overline}
\newcommand{\we}{\wedge}

\newcommand{\FF}{{\cal F}}
\newcommand{\EE}{{\cal E}}

\newcommand{\LL}{{\cal L}}
\newcommand{\MM}{{\cal M}}
\newcommand{\KK}{{\cal K}}

\newcommand{\Hom}{\operatorname{Hom}}
\newcommand{\Ext}{\operatorname{Ext}}

\renewcommand{\a}{\alpha}
\renewcommand{\b}{\beta}
\newcommand{\om}{\omega}
\newcommand{\De}{\Delta}
\newcommand{\la}{\lambda}
\newcommand{\th}{\theta}
\newcommand{\C}{{\Bbb C}}
\newcommand{\R}{{\Bbb R}}
\newcommand{\Z}{{\Bbb Z}}
\newcommand{\Q}{{\Bbb Q}}
\newcommand{\La}{\Lambda}
\newcommand{\Ga}{\Gamma}

\newcommand{\wt}{\widetilde}

\newcommand{\sign}{\operatorname{sign}}

\newcommand{\sub}{\subset}
\newcommand{\ed}{\qed\vspace{3mm}}

\newcommand{\bz}{{\bf z}}

\title{Universal triple Massey products on elliptic curves and 
Hecke's indefinite theta series}
\author{V. Pasol}
\address{Department of Mathematics and Statistics, Boston University,
Boston, MA 02215}
\email{vpasol@@math.bu.edu}
\author{A. Polishchuk}
\address{Department of Mathematics, University of Oregon, Eugene, OR 97403}  
\email{apolish@@math.uoregon.edu}
\thanks{Supported in part by NSF grant}

\begin{document}
\begin{abstract} Generalizing \cite{P1} we express 
universal triple Massey products between line bundles on elliptic curves
in terms of Hecke's indefinite theta series. We show that all Hecke's
indefinite theta series arise in this way.
\end{abstract}
\maketitle

\bigskip

\centerline{\sc Introduction}

\medskip

In general, Massey products are certain operations
on cohomology of dg-algebras (or dg-categories). 
We are interested in triple Massey products 
for morphisms in derived categories of coherent sheaves
on elliptic curves. More precisely, we want to look only at
those Massey products that can be constructed universally and hence give
rise to modular forms. 
A relation between such Massey products and
Hecke's indefinite theta series (introduced and studied by Hecke in \cite{He1}, 
\cite{He2}) was first observed in \cite{P1}. In the present paper
we investigate this connection for a broader class of universal
Massey products.
Our main result is that coefficients of universal Massey products between
line bundles on elliptic curves are always Hecke's indefinite
theta series and that the space
of Hecke's series is spanned by such coefficients (see Theorems \ref{main-thm1} 
and \ref{main-thm2}).

The original definition of Hecke's theta series associated with quadratic forms
of signature $(1,1)$ uses the summation over the cone of positive lattice elements
modulo the action of a subgroup of finite index in the group of automorphisms. In 
this
paper we will use a different set of series described in
\cite{P2}, where the summation is taken over a rational subcone of the
cone of positive lattice elements (see section \ref{def-sec}). 
According to the main theorem of
\cite{P2} these series span the same 
space of $q$-series as Hecke's indefinite theta series.

The main concept that allows to unite Massey products with indefinite theta series
is that of Fukaya product for a configuration of circles in a symplectic torus
(see \cite{F} and \cite{FOOO} for more general discussion of Fukaya categories).
According to the homological mirror conjecture (see \cite{K}), proved for this 
case in 
\cite{PZ} and \cite{P-hmc}, the natural $A_{\infty}$-category structure on
vector bundles over an elliptic curves can be matched with the 
$A_{\infty}$-structure
provided by the Fukaya product. Our Theorem \ref{main-thm1} follows essentially 
from
the observation that the summation pattern over a rank-$2$ lattice appearing in a
Fukaya product corresponding to a quadruple of circles in a symplectic torus,
is the same pattern that is used to form the indefinite theta series of 
\cite{P2}.
To prove Theorem \ref{main-thm2} we have to investigate the
combinatorics of this relation in more detail.

Our results provide a geometric interpretation for all Hecke's
indefinite theta series. We expect that this picture should be
very useful for the study of such series, in the same way as the geometrical
interpretation of usual theta series has become an indispensable tool in their 
study.
In particular, we believe that this framework can be used
to describe all linear relations between indefinite theta series (this problem
was raised by Hecke in \cite{He2}). 

\noindent
{\it Notation.} We use the term $q$-series for formal power series in 
$q^{1/n}$ for
some integer $n>0$. When $q$ is specialized to a complex number $\exp(2\pi 
i\tau)$,
we implicitly assume that $q^{1/n}=\exp(2\pi i\tau/n)$.
For a pair of integer numbers $a$ and $b$ we denote by $gcd(a,b)$ their greatest 
common divisor. We use the shorthand $(\sum_{i\in I_1}-\sum_{i\in I_2})a_i$
for the expression $\sum_{i\in I_1}a_i-\sum_{i\in I_2}a_i$. 

\section{Definitions and results}
\label{def-sec}

\subsection{Hecke's indefinite theta series}\label{indef-theta-sec}

By a {\it lattice} in $\Q^2$ we mean a free abelian subgroup of rank $2$
in $\Q^2$.

\begin{defi} We say that a complex-valued function $f$ on $\Q^2$ is {\it doubly 
periodic}
if there exists a pair of lattices $\La_1\sub\La_2\sub\Q^2$ such that
$f$ is supported on $\La_2$ and is $\La_1$-periodic.
\end{defi}

Let $Q(m,n)=am^2+2bmn+cn^2$ be a $\Q$-valued indefinite quadratic form
on $\Q^2$ (so $b^2>ac$) such that $a$, $b$ and $c$ are positive.
Let $f(m,n)$ be a doubly periodic complex-valued function on $\Q^2$. 
We impose the following condition on $Q$ and $f$:
\begin{equation}\label{f-condition}
f(A\bx)=f(B\bx)=-f(\bx)
\end{equation}
for every $\bx\in\Q^2$, where
$$A=\left(\matrix -1 & -\frac{2b}{a} \\ 0 & 1\endmatrix\right),$$
$$B=\left(\matrix 1 & 0 \\ -\frac{2b}{c} & -1\endmatrix\right).$$
Then we define the 
{\it indefinite theta series associated with} $Q$ and $f$ to be the $q$-series
$$\Th_{Q,f}(q)=(\sum_{m\ge 0,n\ge 0}-
\sum_{m<0,n<0})f(m,n)q^{Q(m,n)/2}.$$
It is easy to see that this series converges for $|q|<1$. 
The equations \eqref{f-condition} are equivalent to the condition that sums of
$f(m,n)q^{Q(m,n)/2}$ over any vertical or horizontal line in $\Q^2$ are zero
(see \cite{P2}). It follows that for every pair of irrational real numbers $c_1$ 
and $c_2$ one has 
$$\Th_{Q,f}(q)=\sum_{(m+c_1)(n+c_2)>0}\sign(m+c_1)f(m,n)q^{Q(m,n)/2}.$$

In \cite{P2} we considered only series $\Th_{Q,f}(q)$ associated with $(Q,f)$ such
that $Q/2$ takes integer values on the support of $f$,
All such series are modular of weight $1$ with respect to some subgroup
$\Ga_0(N)\sub\SL_2(\Z)$. Furthermore, as was shown in \cite{P2} the
space they span coincides with the span of Hecke's indefinite theta series.
In the present work we do note impose the integrality restriction, so we
get series with rational powers of $q$ (but with bounded denominators). Still,
all these series are modular of weight $1$ with respect to some congruence 
subgroup
$\Ga(N)\sub\SL_2(\Z)$.

\subsection{Massey products on a single elliptic curve}\label{Mas-sec}

Let $E$ be an elliptic curve, $L$ be a line bundle of positive
degree $d$ on $E$. Then the space of global sections $H^0(E,L)$ is
an irreducible representation of a certain Heisenberg group $G(L)$ (see \cite{M}
or \cite{M2}).
A {\it theta structure} for $L$ is an isomorphism $G(L)$ with the standard
Heisenberg group $G_d$ generated by $U_1$ and $U_2$ and the central
subgroup $\mu_d=\{\zeta\in\C^*:\ \zeta^d=1\}$ with relations
$$U_1^d=U_2^d=1,\ U_1U_2=\exp(2\pi i/d)U_2U_1.$$
A choice of such a structure defines a canonical (up to rescaling) isomorphism
of $H^0(E,L)$ with the standard model of the irreducible representation of $G_d$.
Hence, we get a canonical basis $(e_0,\ldots,e_{d-1})$ in $H^0(E,L)$ up to a
rescaling $(e_0,\ldots,e_{d-1})\mapsto (\la e_0,\ldots,\la e_{d-1})$.
Namely, $e_0$ is the unique (up to rescaling) vector stabilized by $U_1$ and
$e_k=U_2^k e_0$.

Now let $L_1$, $L_2$ and $L$ be
a triple of line bundles over $E$ of positive degrees $d_1$, $d_2$ and $d$, where  
$d>d_1$, $d>d_2$ and $d_0:=d_1+d_2-d>0$. 
Let us denote $L_0=L_1\otimes L_2\otimes L^{-1}$. 
The triple Massey product associated with this data is a linear map 
$$MP:K\to H^0(L_0),$$ 
where
\begin{align*}
&K=\ker((m_2\otimes\id,\id\otimes m_2):
H^0(L_1)\otimes H^1(L^{-1})\otimes H^0(L_2)\to \\
&H^1(L_1\otimes L^{-1})\otimes
H^0(L_2)\oplus H^0(L_1)\otimes H^1(L^{-1}\otimes L_2)),
\end{align*}
with $m_2$ denoting the natural double product operation.
Let us recall the definition of $MP$ using Dolbeault description of cohomology.
Given an element $\sum_i x_i\otimes y_i\otimes z_i\in K$, we can represent
$y_i$'s by $(0,1)$-forms $\wt{y}_i$ with values in $L^{-1}$. Then
we have
$$\sum_i m_2(x_i,\wt{y}_i)\otimes z_i=\sum_j \ov{\partial}(t_j)\otimes z'_j,$$
$$\sum_i x_i\otimes m_2(\wt{y}_i,z_i)=\sum_k x'_k\otimes\ov{\partial}(s_k)$$
with some $x'_k\in H^0(L_1)$, $z'_j\in H^0(L_2)$,
$t_j\in C^{\infty}(L_1\otimes L^{-1})$ and $s_k\in C^{\infty}(L^{-1}\otimes L_2)$.
Now by definition we have
$$MP(\sum_i x_i\otimes y_i\otimes z_i):=\sum_j m_2(t_j,z'_j)-\sum_k m_2(x'_k,s_k)
\in H^0(L_0).$$
One can also give a purely algebraic definition (see section \ref{Mas-cont-sec}).

Let us assume in addition that all the line bundles $L_1$, $L_2$, $L$ and $L_0$
are equipped with theta structures. Then we can choose some
canonical bases $(e_1,\ldots,e_{d_1})$ in $H^0(E,L_1)$, $(f_1,\ldots,f_{d_2})$ in 
$H^0(E,L_2)$, $(g_1,\ldots,g_d)$ in $H^1(E,L^{-1})$, and 
$(h^1,\ldots,h_{d_0})$ in $H^0(E,L_0)$.

\begin{defi} A collection of complex constants $\bc=(c_{ijk})$, where
$(i,j,k)$ varies through $\Z/d_1\Z\times\Z/d\Z\times\Z/d_2\Z$, is called
a {\it Massey system} for $(L_1,L_2,L)$
if we have $\sum_{ijk}c_{ijk}e_i\otimes f_j\otimes g_k\in K$.
This is equivalent to the equation
\begin{equation}\label{Mas-van-cond}
(m_2\otimes\id,\id\otimes m_2)(\sum_{i,j,k}c_{ijk}e_i\otimes g_j\otimes f_k)=0.
\end{equation}
\end{defi}

Note that the condition \eqref{Mas-van-cond} does not depend on the choices of 
canonical bases. For every Massey system $\bc$ the corresponding Massey product 
\begin{equation}\label{Mas-element}
MP(\sum_{i,j,k}c_{ijk}e_i\otimes g_j\otimes f_k)\in H^0(L_0)
\end{equation}
is well defined and hence we can also consider coefficients
\begin{equation}\label{Mas-coef}
\lan MP(\sum_{i,j,k}c_{ijk}e_i\otimes g_j\otimes f_k),h_l\ran, \ l\in\Z/d_0\Z
\end{equation}
of this element with respect to the base $(h_l)$.
Below we are going to study the relative analogue of this construction when
the elliptic curve and the line bundles on it vary in a family.

\subsection{Universal Massey products}\label{univ-Mas-sec}

Now let $\pi:\EE\to S$ be a family of elliptic curves and $\LL$ 
a line bundle on $\EE$ of relative degree $d$. 
Assume that a relative theta structure for $\LL$ is chosen.
Then we have
a $\G_m$-torsor $\Bases(\LL)$ over $S$ such that its fiber of $s\in S$
is the set of canonical bases in $H^0(E_s,L_s)$, where $E_s\sub\EE$
is the elliptic curve corresponding to $s\in S$ and $L_s=\LL|_{E_s}$.
Let $L\Bases(\LL)$ be the line bundle over $S$ associated with 
$\G_m$. Note that we also have $L\Bases(\LL)\simeq (\pi_*L)^{U_1}$,
where the action of $U_1$ on $\pi_*L$ is induced by the theta structure for $\LL$.
The degree $d$ map $\Bases(\LL)\to\det\pi_*\LL$ sending a basis
$e_1,\ldots,e_d$ to $e_1\we\ldots\we e_d$ induces an isomorphism
\begin{equation}\label{det-isom}
L\Bases(\LL)^d\simeq\det R\pi_*\LL,
\end{equation}
where $\det R\pi_*\LL$ is the determinant line bundle on $S$
(since $d>0$, it is the determinant of the vector bundle $\pi_*\LL$).
Similarly, if $\LL$ is a line bundle on $\EE$ of negative relative degree $d$
and $L\Bases(\LL)$ is the line bundle on $S$ associated with the $\G_m$-torsor
of canonical bases in $H^1(E_s,L_s)$, then the isomorphism \eqref{det-isom}
still holds.

Now assume that we have a triple of line bundles $\LL_1$, $\LL_2$ and
$\LL$ over $\EE$ of positive relative degrees $d_1$, $d_2$ and $d$, where  
$d>d_1$, $d>d_2$ and $d_0:=d_1+d_2-d>0$. 
Let us denote
$\LL_0=\LL_1\otimes\LL_2\otimes\LL^{-1}$. Assume that all the line bundles
$\LL_1$, $\LL_2$, $\LL$ and $\LL_0$ are equipped with relative theta structures.
Then for every point $s\in S$ we can consider the triple
Massey product on $E_s$ associated with the restrictions of 
$\LL_1$, $\LL_2$ and $\LL$ to $E_s$.
Furthermore, as we will show later (see Lemma \ref{Mas-cont-lem}), there exists a 
vector bundle $\KK$ over $S$ such that $\KK|_s$ is the
source of the corresponding Massey product, so that we have a global morphism
of vector bundles
$$MP:\KK\to \pi_*\LL_0$$ 
inducing the Massey product at every point.
Now we want to use relative theta structures on our bundles to produce
line subbundles in $\KK$. Recall that these structures allow us to
choose canonical bases of various $H^0$ and $H^1$ spaces over members
$E_s$ of our family. 

\begin{defi}
We say that a collection of constants $\bc=(c_{ijk})$ (independent of $s\in S$), 
where
$(i,j,k)$ varies through $\Z/d_1\Z\times\Z/d\Z\times\Z/d_2\Z$,
is a {\it universal Massey system for} $(\LL_1,\LL_2,\LL)$ if $\bc$ is a
Massey system for the restrictions of these line bundles to $E_s$ for every $s$.
Equivalently,
the vanishing condition \eqref{Mas-van-cond} should hold for all $s\in S$. For 
every
such a system the elements \eqref{Mas-element} define a morphism
$$MP(\bc):L\Bases(\LL_1)\otimes L\Bases(\LL^{-1})\otimes L\Bases(\LL_2)\to 
\pi_*\LL_0$$
that we call a {\it universal Massey product over $S$}. Similarly, constants
\eqref{Mas-coef} define sections of a certain line bundle on $S$:
$$MP(\bc)_l\in H^0(S,L\Bases(\LL_1)^{-1}\otimes L\Bases(\LL^{-1})^{-1}\otimes
L\Bases(\LL_2)^{-1}\otimes L\Bases(\LL_0))$$
that we call {\it coefficients} of the above universal Massey product over $S$.
\end{defi}

Let $e:S\to \EE$ be the relative neutral point on $\EE$ and let us set
$\ov{\om}=e^*\om_{\EE/S}$, where $\om_{\EE/S}$ is the relative canonical bundle.
We assume in addition that all our line bundles $\LL_1$, $\LL_2$ and $\LL$
are trivialized along the zero section and are symmetric up to torsion, i.e.,
$[-1]^*\LL\otimes\LL^{-1}$ is torsion, etc.

\begin{lem} In the above situation one has the following equality in
the group $\Pic(S)\otimes\Q$
$$-[L\Bases(\LL_1)]-[L\Bases(\LL_2)]-[L\Bases(\LL^{-
1})]+[L\Bases(\LL_0)]=[\ov{\om}].$$
\end{lem}

\Pf . In $\Pic(S)\otimes\Q$ one has 
$$[\det R\pi_*\LL_1]=-d_1\cdot[\ov{\om}]/2$$
(see e.g., \cite{MB}, ch.VIII).
Hence, $[L\Bases(\LL_1)]=-[\ov{\om}]/2$. Applying similar formulas for
other line bundles we immediately get the result.
\ed

The above computation makes it reasonable to conjecture that
coefficients of universal Massey products for appropriate line bundles
$\LL_1$, $\LL_2$ and $\LL$ over some moduli stack of elliptic curves (with 
additional structures) define modular forms of
weight $1$ with respect to appropriate congruence subgroups. 
Furthermore, it is plausible that the techniques of \cite{MB} can be used to
prove this algebraically. However, we will check this by directly expressing
these coefficients in terms of Hecke's indefinite theta series.

\subsection{Formulation of the results}\label{results-sec}

Let $\EE$ be the standard family of elliptic curves parametrized by the
upper half-plane $\gH$, i.e., $\EE=(\gH\times\C)/\Z^2$,
where $(m,n)\in\Z^2$ acts by $(\tau,z)\mapsto (\tau,z+m\tau+n)$. 
There is a natural action of
$\SL_2(\Z)$ on the family $\EE\to\gH$ given by
$$g(\tau,z)=(\frac{a\tau+b}{c\tau+d},\frac{z}{c\tau+d}),$$
where $g=\left(\matrix a & b \\ c & d \endmatrix\right)$.
We are going to consider an equivariant
version of the previous construction so that the resulting universal Massey 
product
will be a section of a line bundle on the quotient-stack of $\gH$ by the
an appropriate congruenz-subgroup of $\SL_2(\Z)$. 

For every integer $N$ and a pair of real numbers $(r,s)$ we
denote by $\LL(N)[r,s]$ the line bundle on $\EE$ obtained by taking the
quotient of the trivial bundle $\gH\times\C\times\C$ over $\gH\times\C$
by the following action of $\Z^2$:
$$(m,n)\cdot(\tau,z,\la)=(\tau,z+m\tau+n,\la\cdot
\exp(-\pi iNm^2\tau-2\pi iNmz-2\pi im(r\tau+s))).$$
When we fix $\tau\in\gH$ we will denote by $L(N)[r,s]$ the line bundle of degree 
$d$ on the elliptic curve $E_{\tau}=\C/(\Z\tau+\Z)$ obtained from $\LL(N)[r,s]$ by 
restriction. For the remainder of this section we assume that $r$ and $s$ are 
rational. 
Then the line bundle $\LL(N)[r,s]$ has a natural equivariant structure with
respect to the subgroup $\Ga_{1,2}\cap\Ga[r,s]\sub\SL_2(\Z)$, where
$\Ga_{1,2}\sub\SL_2(\Z)$ consists of matrices such that $ab$ and $cd$ are
even, and $\Ga[r,s]$ is the subgroup consisting of $g\in\SL_2(\Z)$ such that
$$g^t\left(\matrix r \\ s\endmatrix\right)-\left(\matrix r \\ 
s\endmatrix\right)\in\Z^2.$$
This equivariant structure corresponds to the following action of
the above subgroup on $\gH\times\C\times\C$:
$$g(\tau,z,\la)=(\frac{a\tau+b}{c\tau+d},\frac{z}{c\tau+d},\la\cdot
\exp[\pi iN\frac{cz^2}{c\tau+d}+2\pi irz(1-\frac{1}{c\tau+d})]).$$
It is easy to see that we have canonical isomorphisms
$$\LL(N+N')[r+r',s+s']\simeq\LL(N)[r,s]\otimes\LL(N')[r',s']$$
compatible with equivariant structures.

Note that the line bundle $\LL(N)[r,s]$ has relative degree $N$.
For $N\neq 0$ the natural theta structure on it is given by the following
action of the standard Heisenberg group $G_N$:
$$U_1(\tau,z,\la)=(\tau,z-1/N,\la),$$
$$U_2(\tau,z,\la)=(\tau,z-\tau/N,\la\cdot\exp(-\pi i\frac{\tau}{N}+2\pi i z+
2\pi i\frac{r\tau+s}{N})).$$
One can check its compatibility with the action of an appropriate  
congruenz-subgroup of
$\SL_2(\Z)$.
The corresponding action of $G_N$ on sections of $\LL(N)[r,s]$
is given by 
$$U_1f(\tau,z)=f(\tau,z+1/N)$$
$$U_2f(\tau,z)=f(\tau,z+\frac{\tau}{N})\exp(\pi i\frac{\tau}{N}+2\pi i z+
2\pi i\frac{r\tau+s}{N}).$$
One can easily check that for $N>0$ the functions
$$\th_{N,k}[r,s](\tau,z)=\sum_{m\in N\Z+k}\exp(\pi i\frac{\tau}{N}m^2+2\pi im(z+
\frac{r\tau+s}{N})),$$
where $k=0,\ldots,N-1$, 
descend to a canonical basis of the space of global sections
of the restriction of 
$\LL(N)[r,s]$ to every member of the family $\EE\to\gH$.
Furthermore, the functional equation for theta functions shows that the functions
\begin{equation}\label{basis-eq}
F_{N,k}[r,s]:=\exp(\pi i\frac{\tau}{N}r^2)\cdot\th_{N,k}[r,s], \ k\in\Z/N\Z,
\end{equation}
descend to sections of 
$\LL(N)[r,s]\otimes \MM$,  
where $\MM$ is a line bundle on $\gH$ equipped with an action of
some congruence subgroup of $\SL_2(\Z)$ and an equivariant morphism
$\MM^{8}\simeq\om_{\EE/\gH}^2$. 
Thus, we obtain
an equivariant isomorphism of line bundles 
\begin{equation}\label{F-isom}
L\Bases(\LL(N)[r,s])\simeq\ov{\om}^{-1/2}
\end{equation}
By Serre duality we also get similar isomorphism for $N<0$. 

Let $(d_1,d_2,d)$ be positive integers satisfying assumptions
of section \ref{univ-Mas-sec}. Let us also fix a quadruple of rational numbers
$(v_1,v_2,w_1,w_2)$. With these data we associate line bundles 
$\LL_1$, $\LL_2$ and $\LL_0$ on $\EE$ as follows:
\begin{equation}\label{line-bun-eq}
\begin{array}{l}
\LL_1=\LL(d_1)[d_1v_1,d_1w_1],\\
\LL_2=\LL(d_2)[-d_2v_2,-d_2w_2],\\
\LL_0=\LL(d_0)[0,0],
\end{array}
\end{equation}
We set $\LL=\LL_1\otimes\LL_2\otimes\LL_0^{-1}$, so that there
is a canonical isomorphism
$$\LL\simeq\LL(d)[d_1v_1-d_2v_2,d_1w_1-d_2w_2].$$
Using isomorphisms \eqref{F-isom}
coefficients of universal Massey products associated with $(\LL_1,\LL_2,\LL)$
can be viewed as functions on the upper half-plane satisfying the modular 
functional
equation of weight $1$. The following theorem identifies these coefficients with 
certain
Hecke's indefinite theta series.

\begin{thm}\label{main-thm1} 
For every universal Massey system $\bc$ for $(\LL_1,\LL_2,\LL)$
and every $l\in\Z/d_0\Z$ we have
$$MP(\bc)_l=\Th_{Q,f_{\bc,l}}(q),$$
where $q=\exp(2\pi i\tau)$,
$$Q(m,n)=Q_{d_1,d_2,d}(m,n)=\frac{1}{d}[d_1(d-d_1)m^2+2d_1d_2mn+d_2(d-d_2)n^2],$$
and $f_{\bc,l}$ is a certain doubly periodic function on $\Q^2$ such that the
data $(Q,f_{\bc,l})$ satisfy the conditions of section \ref{indef-theta-sec}.
\end{thm}

The explicit form of the function $f_{\bc,l}$ is given by \eqref{f-eq}.
Our second result shows that every Hecke's series comes from a
universal Massey product.

\begin{thm}\label{main-thm2}
For every data $(Q,f)$ as in \ref{indef-theta-sec} there exist positive
integers $(d,d_1,d_2)$ satisfying the inequalities
$$d-d_1>0,\ d-d_2>0,\ d_1+d_2-d>0, $$
and a universal Massey system $\bc$ for $(\LL_1,\LL_2,\LL)$, where 
$$\LL_1=\LL(d_1)[0,0],\ \LL_2=\LL(-d_2)[0,0],\ \LL=\LL(d)[0,0],$$
such that
$$\Th_{Q,f}(q)=\sum_{l\in\Z/d_0\Z} MP(\bc)_l.$$
\end{thm}

\subsection{Continuity of Massey products}
\label{Mas-cont-sec}

In this section we recall an algebraic definition of Massey products considered 
above. As a consequence we will derive that they vary continuously with 
parameters.

First, we observe that in the situation of section \ref{Mas-sec} there is an 
isomorphism
$$H^0(L_1)\otimes H^1(L^{-1})\otimes H^0(L_2)\simeq\Ext^1(X,Y),$$ 
where 
$X=H^0(L_1)^*\otimes L_1$ and $Y=H^0(L_2)\otimes L_1\otimes L^{-1}$.
Let
$$0\to Y\to V\to \Ext^1(X,Y)\otimes X\to 0$$
be the universal extension, and let
$$0\to Y\to V_K\to K\otimes X\to 0$$ be the  induced
extension corresponding to the inclusion $K\sub\Ext^1(X,Y)$.
Note that we have natural morphisms $\a:\OO_E\to X$ and 
$\b:Y\to L_2\otimes L_1\otimes L^{-1}\simeq L_0$.
It is easy to check that $K$ coincides with the subspace of elements 
$e\in\Ext^1(X,Y)$ such that $e\circ\a=0$ and $\b\circ e=0$. Furthermore, one has
$$MP(e)=MP(\a,e,\b),$$
where the Massey products in the RHS are of the type considered in 
\cite{P1}, sec. 6.1. Hence, from Proposition 6.1 of loc.~cit. we obtain the
following interpretation of the map $MP$. Let 
$$\wt{\a}:K\otimes\OO_E\to V_K$$ 
be the unique lifting of the map
$\id\otimes\a:K\otimes\OO_E\to K\otimes X$, and let
$$\wt{\b}:V_K\to L_0$$ 
be the unique map extending $\b:Y\to L_0$.
Then the Massey product $MP:K\to H^0(L_0)$ coincides with the map of the spaces of 
global sections induced by the morphism
$$\wt{\b}\circ\wt{\a}:K\otimes\OO_E\to L_0.$$ 

Now assume that we are in the situation of section \ref{univ-Mas-sec}, so we have 
a family $\pi:\EE\to S$ of elliptic curves equipped with line bundles
$(\LL_1,\LL_2,\LL)$. Then the
relative version of the above considerations leads to the following result.

\begin{lem}\label{Mas-cont-lem} 
There exists a vector bundle $\KK$ over $S$ and a morphism
of vector bundles
$MP:\KK\to \pi_*\LL_0$ which restricts to the above morphism $MP:K\to H^0(E,L_0)$
on every fiber of $\pi$.
\end{lem}

\Pf . Our inequalities on $d_1$, $d_2$ and $d$ imply that $\max(d_1,d_2)>1$.
Without loss of generality we can assume that $d_2>1$.
Now we are going to rewrite the definition of $K\sub\Ext^1(X,Y)$ for a single
elliptic curve $E=\pi^{-1}(s)$ in such a form that it will become clear that it 
produces a vector
bundle $\KK$ over $S$.
Let us define coherent sheaves $X'$ and $Y'$ on $E$ from exact sequences
$$0\to\OO_E\stackrel{\a}{\to} X\to X'\to 0,$$
$$0\to Y'\to Y\stackrel{\b}{\to} L_0\to 0$$
(surjectivity of $\b$ follows from the assumption $d_2>1$).
Since $d_1>d_0$, we have $\Hom(X,L_0)=0$, hence the kernel of the map
$$\Ext^1(X,Y)\to\Ext^1(X,L_0):e\mapsto \b\circ e$$
can be identified with $\Ext^1(X,Y')$. Now by definition, $K$ is the kernel
of the map
\begin{equation}\label{K-map}
\Ext^1(X,Y')\to\Ext^1(\OO_E,Y) 
\end{equation}
induced
by $\a$ and by the embedding $Y'\to Y$.
Since $d_1<d$, one has $\Hom(\OO_E,Y')\sub\Hom(\OO_E,Y)=0$. 
Hence, we have exact sequences
$$0\to\Ext^1(X',Y')\to\Ext^1(X,Y')\stackrel{g}{\to}\Ext^1(\OO_E,Y')\to  0,$$
$$0\to\Hom(\OO_E,L_0)\to\Ext^1(\OO_E,Y')\stackrel{f}{\to}\Ext^1(\OO_E,Y)\to  0.$$
The map \eqref{K-map} is equal to the composition $f\circ g$, so from these exact 
sequences we get the following description of $K$ as an extension:
$$0\to\Ext^1(X',Y')\to K\to H^0(L_0)\to 0.$$
Since $\Hom(X',Y')\sub\Hom(X,Y)=0$, we see that there is a vector bundle $\KK'$ 
over $S$ with the fiber $\Ext^1(X',Y')$ over $s$. Now $\KK$ can be constructed as 
an extension of $\pi_*\LL_0$ by $\KK'$.
\ed

Later we will apply this lemma for the family of line bundles
\eqref{line-bun-eq} depending on real parameters $(v_1,v_2)$
(keeping the elliptic curve and the parameters $(w_1,w_2)$ constant). Then the
corresponding Massey systems vary in a vector bundle over $\R^2$ and
the Massey product is a continuous map from the total space of this bundle
to the vector space $H^0(L_0)$.

\section{Calculations}

We follow the method developed in \cite{P1}.
The basic idea is to use homological mirror symmetry to relate the Massey products
considered above to triple Fukaya products for symplectic tori, so we start by
reviewing this relation.

\subsection{Connection with Fukaya products}

Let $\tau$ be an element of the upper half-plane and let
$E=E_{\tau}$ be the corresponding elliptic curve.
The mirror partner of $E$ is the torus $T=\R^2/\Z^2$ with the complexified
symplectic form $\om=-2\pi i\tau dx\we dy$. 
We will work with the subcategory $\FF_s$ in the Fukaya category of $(T,\om)$ 
described in section 2.1 of \cite{P1}. Recall that objects of
$\FF_s$ are pairs $(\ov{\ell},t)$, where $\ov{\ell}\sub T$ is the image of a 
nonvertical line $\ell\in\R^2$ of rational slope, and $t$ is a real number. The 
space of morphisms between $(\ov{\ell_1},t_1)$ and $(\ov{\ell_2},t_2)$ is defined 
only when
$\ov{\ell_1}\neq\ov{\ell_2}$ and is set to be 
$\oplus_{P\in\ov{\ell_1}\cap\ov{\ell_2}}\C[P]$. This space has degree $0$ if
the slope of $\ell_1$ is smaller than the slope of $\ell_2$ and degree $1$ 
otherwise. By definition $m_1=0$ while
the operation $m_k$ for $k\ge 2$ is defined using certain summation over 
(k+1)-gons (see \cite{P1}, sec. 2.1). Below we will only 
use operations $m_2$ (preserving degree) and $m_3$ (lowering degree by 
$1$). A little subtlety here is that operations $m_k$ are defined only for
transversal configurations, so one has to change the axiomatics slightly
(see \cite{KS}, sec. 4.3).

The main theorem of \cite{P-hmc} extending the result of \cite{PZ} gives an 
$A_{\infty}$-equivalence between $\FF_s$ and the natural $A_{\infty}$-category 
whose objects are stable bundles on $E$ (the latter $A_{\infty}$-structure is 
obtained from the dg-category structure given by the standard enhancement on
the derived category of coherent sheaves). Let us recall how this equivalence is 
defined 
for line bundles on $E$ (which correspond to lines of integer slope on $T$).

A line bundle $L(N)[r,s]$ on $E$ (see section \ref{results-sec}) corresponds to a 
pair $(\ov{\ell_{N,r}},-s)$, where $\ell_{N,r}$ is the
line $\{(x,Nx-r), x\in\R\}\sub\R^2$.
The identification of morphism spaces is defined as follows.
Assume that $N_1<N_2$. Then we have
$$\Hom(L(N_1)[r_1,s_1],L(N_2)[r_2,s_2])\simeq  H^0(E,L(N_2-N_1)[r_2-r_1,s_2-
s_1]),$$
$$\Hom_{\FF_s}((\ov{\ell_{N_1,r_1}},-s_1),(\ov{\ell_{N_2,r_2}},-s_2)=
\oplus_{k\in\Z/(N_2-N_1)\Z}\C[P_k],$$
where 
$$P_k=(\frac{k+r_2-r_1}{N_2-N_1},\frac{N_1k+N_1r_2-N_2r_1}{N_2-N_1})\in
\ov{\ell_{N_1,r_1}}\cap\ov{\ell_{N_2,r_2}}$$
are points of intersection between the corresponding lines on $T$.
The isomorphism between the above two morphism spaces sends the basis of theta 
functions
$$\th_{N_2-N_1,k}[r_2-r_1,s_2-s_1],\ k\in\Z/(N_2-N_1)\Z,$$
of the former space to the basis
\begin{equation}\label{Fuk-basis-eq}
e(k)=e_{N_1,N_2}[r_1,r_2;s_1,s_2](k):=\exp(-\pi i \frac{\tau(r_2-r_1)^2}{N_2-N_1}-
2\pi i \frac{(s_2-s_1)(r_2-r_1)}{N_2-N_1})[P_k]
\end{equation}
of the latter space.

In the case $N_2>N_1$ the spaces of morphisms have degree $1$ and the 
identification uses Serre duality to reduce to the previous case.

Our Massey product coincides with the restriction of the triple product
$$m_3:H^0(E,L_1)\otimes H^1(E,L^{-1})\otimes H^0(E,L_2)\to H^0(L_0)$$
to $K$ (see \cite{P3}, sec. 1.1).
Also, Massey products are preserved by $A_{\infty}$-equivalences
(see e.g. Proposition 1.1 of \cite{P3}).
Hence, using the above equivalence with $\FF_s$ we derive that in the situation of 
section \ref{results-sec} the coefficients 
$MP(\bc)_l(\tau)$ can be computed as follows. Let $L$ (resp., $L_0$, $L_1$, $L_2$) 
be the restriction of $\LL$ (resp., $\LL_0$, $\LL_1$, $\LL_2$) to $E=E_{\tau}$.
Then our Massey product can be expressed in terms of the triple product
$$m_3:\Hom_{\FF_s}(O_0,O_1)\otimes\Hom_{\FF_s}(O_1,O_2)\otimes
\Hom_{\FF_s}(O_2,O_3)\to\Hom_{\FF_s}(O_0,O_3),$$
where $(O_0,O_1,O_2,O_3)$ is the quadruple of objects in the Fukaya category 
corresponding to the line bundles $(\OO_E,L_1,L_1\otimes L^{-1}\simeq 
L_2^{-1}\otimes L_0, L_0)$.
More precisely, 
using the above dictionary we find
$$O_0=(\ov{\ell_{0,0}},0),$$
$$O_1=(\ov{\ell_{d_1, d_1v_1}},-d_1w_1),$$
$$O_2=(\ov{\ell_{d_1-d,d_2v_2}},-d_2w_2),$$
$$O_3=(\ov{\ell_{d_0,0}},0).$$
Assume that $v_1$ and $v_2$ are sufficiently generic, so that
these $4$ objects form a transversal configuration on the torus
$\R^2/Z^2$ (i.e., no three of them have a common intersection 
point).  Consider the bases $(e_{ab}(k))$ 
of the spaces $\Hom_{\FF_s}(O_a,O_b)$ defined by \eqref{Fuk-basis-eq}. Then
the basis \eqref{basis-eq} of $H^0(E,L_1)$
corresponds to 
$$\exp(\pi i\tau d_1v_1^2)\cdot e_{01}(i), \ i\in\Z/d_1\Z.$$
Similarly finding the correct bases of other spaces we arrive to the following

\begin{lem}\label{connection-lem}
A collection $\bc=(c_{ijk})$ is a Massey system for $(L_1,L_2,L)$ iff 
the system of equations
\begin{equation}
\sum_{ij}c_{ijk}m_2(e_{01}(i),e_{12}(j))=0, \ k\in\Z/d_2\Z,
\end{equation}
\begin{equation}
\sum_{jk}c_{ijk}m_2(e_{12}(j),e_{23}(k))=0, \ i\in\Z/d_1\Z
\end{equation}
holds.
Assume that $v_1$ and $v_2$ are generic. Then for a Massey system $\bc$ one has
$$MP(\bc)_l(\tau)=\exp(\pi i\tau Q(v_1,v_2))\cdot
\sum_{ijk}c_{ijk}\lan m_3(e_{01}(i),e_{12}(j),e_{23}(k)),
e_{03}(l)\ran,$$
where $Q$ is the quadratic form introduced in Theorem \ref{main-thm1}.
\end{lem}

Now our plan is to compute explicitly the double and triple products appearing
in this lemma. After rewriting the above formula for $MP(\bc)_l(\tau)$
in a suitable form we will be able to get rid of the genericity assumption on
$v_1$ and $v_2$. Finally, we will apply these considerations to the universal
Massey systems.

\subsection{Double products}\label{double-sec}

From now on we are going to switch to
a slightly different notation for theta functions (similar to the one adopted in 
\cite{P1}): for a subset $I\sub\Q$ of the form $I=a\Z+b$, where $a,b\in\Q$
and $a\neq 0$, we set
$$\th_I(z,\tau)=\sum_{m\in I}\exp(\pi i\tau m^2+2\pi imz) .$$
The relation with the notation of section \ref{results-sec} is the following:
$$\th_{N,k}[r,s](z,\tau)=\th_{N\Z+k}(z+\frac{r\tau+s}{N},\frac{\tau}{N}).$$

Let us set $p_1=dd_1/(d-d_1)$ and $p_2=dd_2/(d-d_2)$.

\begin{lem} For $i\in\Z/d_1\Z$, $j\in\Z/d\Z$, and $k\in\Z/d_2\Z$ 
one has
$$m_2(e_{01}(i),e_{12}(j))=\sum_{u\in\Z/(d-d_1)\Z}
\th_{I_1(u,i,j)}(p_1x_1,p_1\tau)e_{02}(u),$$
$$m_2(e_{12}(j),e_{23}(k))=\sum_{v\in\Z/(d-d_2)\Z}
\th_{I_2(v,j,k)}(p_2x_2,p_2\tau)e_{13}(v),$$
where 
$$x_1=-\frac{[(d-d_1)v_1+d_2v_2]\tau+(d-d_1)w_1+d_2w_2}{d}, $$
$$x_2=-\frac{[(d-d_2)v_2+d_1v_1]\tau+(d-d_2)w_2+d_1w_1}{d},$$ 
$$I_1(u,i,j)=\{m\in\Q:\ m\equiv -\frac{i}{d_1}-\frac{j}{d}\mod(\Z),\ 
\frac{dm}{d-d_1}\equiv -\frac{i}{d_1}-\frac{u}{d-d_1}\mod(\Z)\}$$
$$I_2(v,j,k)=\{m\in\Q:\ m\equiv \frac{j}{d}+\frac{k}{d_2}\mod(\Z),\ 
\frac{dm}{d-d_2}\equiv \frac{v}{d-d_2}+\frac{k}{d_2}\mod(\Z)\}.$$
\end{lem}

\Pf . Formula (2.5) of \cite{P1} in our case gives
$$m_2(e_{01}(i),e_{12}(j))=\sum_{m\in\Z/I_1\Z}
\th_{I_1,-\frac{i}{d_1}-\frac{j}{d}+m}(p_1x_1,p_1\tau)e_{02}(i+j-md),$$
where $I_1=I_1(0,0,0)$.
This easily implies the first identity. The second formula is derived in the 
same way.
\ed

In the remainder of this section we assume that $v_1=v_2=w_1=w_2=0$.
Then a collection $(c_{ijk})$ is a universal Massey system iff the following
system of equations is satisfied identically in $\tau$:
\begin{equation}\label{double-eq1}
\sum_{ij}c_{ijk}\th_{I_1(u,i,j)}(0,p_1\tau)=0 \text{ where }k\in\Z/d_2\Z,
u\in\Z/(d-d_1)\Z,
\end{equation}
\begin{equation}\label{double-eq2}
\sum_{jk}c_{ijk}\th_{I_2(v,j,k)}(0,p_2\tau)=0 \text{ where }i\in\Z/d_1\Z,
v\in\Z/(d-d_2)\Z.
\end{equation}
We would like to rewrite these equations in the form similar to 
\eqref{f-condition}.
Set $\bd\Z^3:=d_1\Z\times d\Z\times d_1\Z$. We can view $\bc=(c_{ijk})$ as a 
$\bd\Z^3$-periodic function on $\Q^3$ supported on
$\Z^3$, by setting $c_{ijk}=0$ for $(i,j,k)\in\Q^3\setminus\Z^3$.
Then the coefficient with $q^{p_1m^2/2}$ in the LHS of \eqref{double-eq1} 
comes from the terms in the sums defining $\th_{I_1(u,i,j)}$ corresponding
to $m$ and to $-m$. It is easy to check that $I_1(u,i,j)$ is nonempty iff
$u\equiv i+j\mod(g_1\Z)$, where $g_1=gcd(d_1,d)$.
Note also that if we fix $m$ and $u$
then the condition $m\in I_1(u,i,j)$
determines the pair $(i,j)\in\Q/d_1\Z\times\Q/d\Z$ uniquely. 
An easy computation shows that if $m\in I_1(u,i,j)$ then
$$-m\in I_1(u,\frac{-2d_1u}{d-d_1}-i,\frac{2du}{d-d_1}-j).$$
Hence, \eqref{double-eq1} is equivalent to the following system of equations
on a $\bd\Z^3$-periodic function $\bc$ on $\Q^3$ supported on
$\Z^3$:
\begin{equation}\label{new-double-eq1}
c_{i,j,k}=-c_{-\frac{2d_1u}{d-d_1}-i,\frac{2du}{d-d_1}-j,k},
\end{equation}
for every $(i,j,k)\in\Q^3$ and $u\in\Z/(d-d_1)\Z$ such that 
$i+j\equiv u\mod(g_1\Z)$. 
Comparing the above equations for $(i,j,k,u)$ and $(i,j,k,g_1+u)$ we derive
that $\bc$ should satisfy the following additional periodicity:
\begin{equation}\label{per-eq1}
c_{i,j,k}=c_{i-\frac{2d_1g_1}{d-d_1},j+\frac{2dg_1}{d-d_1},k}.
\end{equation}
Since $\bc$ is supported on $\Z^3$ this is possible only if
$\frac{2d_1g_1}{d-d_1}$ is an integer (then $\frac{2dg_1}{d-d_1}$ is
necessarily also an integer).
Once the periodicity \eqref{per-eq1} is satisfied we can replace $u$ with
$i+j$ in \eqref{new-double-eq1} and get an equivalent condition
\begin{equation}\label{vanish-oper-eq1}
c_{i,j,k}=
-c_{-\frac{d+d_1}{d-d_1}i-\frac{2d_1}{d-d_1}j,\frac{2d}{d-d_1}i+\frac{d+d_1}{d-
d_1}j,k}
\end{equation}
for $(i,j,k)\in\Q^3$.
Repeating this procedure with equation \eqref{double-eq2} we obtain the 
periodicity
condition
\begin{equation}\label{per-eq2}
c_{i,j,k}=c_{i,j+\frac{2dg_2}{d-d_2},k-\frac{2d_2g_2}{d-d_2}},
\end{equation}
where $g_2=gcd(d_2,d)$, together with the equation
\begin{equation}\label{vanish-oper-eq2}
c_{i,j,k}=
-c_{i,\frac{d+d_2}{d-d_2}j+\frac{2d}{d-d_2}k,-\frac{2d_2}{d-d_2}j-\frac{d+d_2}{d-
d_2}k}.
\end{equation}
Summarizing we obtain the following result.

\begin{lem}\label{tech-lem}
Assume that $v_1=v_2=w_1=w_2=0$. Then 
for a $\bd\Z^3$-periodic
function $\bc$ on $\Q^3$ supported on $\Z^3$
the system of equations 
\eqref{double-eq1}, \eqref{double-eq2} 
is equivalent to the system
\eqref{per-eq1}, \eqref{vanish-oper-eq1}, \eqref{per-eq2}, 
\eqref{vanish-oper-eq2}.
A nonzero solution exists only if
$$\frac{2d_1g_1}{d-d_1}\in\Z \text{ and } \frac{2d_2g_2}{d-d_2}\in\Z,$$
where $g_1=gcd(d_1,d)$ and $g_2=gcd(d_2,d)$.
\end{lem}

\subsection{Triple product}\label{triple-sec}

Define the lattice 
$$\La=\{(\la_1,\la_2)\in\Q^2\ |\ d_1\la_1+(d-d_2)\la_2\in d\Z\text{ and }
(d-d_1)\la_1+d_2\la_2\in d\Z\}$$
and the sublattice $\La^0=\La\cap\Z^2$. 
Let us also consider 
the cone $C=\{(x_1,x_2)\in\R^2:\ x_1x_2>0\}$.

\begin{lem}\label{triple-prod-lem} Assume that $v_1$ and $v_2$ are generic.
Then for a collection of constants $\bc=(c_{ijk})$
one has
$$\sum_{ijk}c_{ijk}m_3(e_{01}(i),e_{12}(j),e_{23}(k))=\sum_{i,j,k,l\in\Z/d_0\Z}
\Th_{ijkl}e_{03}(l)$$ 
with
$$\Th_{ijkl}=
\sum_{\bla\in\La^0(i,j,k,l)\cap (C-\bv)}
\sign(\la_1+v_1)\exp[\pi i\tau Q(\bla)+2\pi i\bla\cdot(\tau\bv+\bw)],$$
where $Q$ is the quadratic form appearing in Theorem \ref{main-thm1},
$\bx\cdot \by=1/2(Q(\bx+\by)-Q(\bx)-Q(\by))$ is the associated bilinear form,
$\bv=(v_1,v_2)$, $\bw=(w_1,w_2)$, and
\begin{align*}
\La^0(i,j,k,l)=\{(\la_1,\la_2)\in\Z^2:\ 
&\la_1\equiv\frac{i}{d_1}-\frac{l}{d_0}\mod(\Z),\\
&\la_2\equiv\frac{l}{d_0}-\frac{k}{d_2}\mod(\Z),\\
&\frac{d_1\la_1+(d-d_2)\la_2}{d}\equiv-\frac{j}{d}-\frac{k}{d_2}\mod(\Z)\}
\end{align*}
\end{lem}

\Pf . Applying formula (2.9) from \cite{P1} we find
$$\sum_{i,j,k}c_{ijk}m_3(e_{01}(i),e_{12}(j),e_{23}(k))=
\sum_{(m,n)\in\La/\La^0}\Th_{mn}
e_{03}(i+j+k+d_0n),$$
where
$$\Th_{mn}=\sum_{\bla\in(\La^0 +(m,n)+\bu(i,j,k))\cap (C-\bv)}\sign(\la_1+v_1)
\exp[\pi i\tau Q(\bla)+2\pi i\bla\cdot(\tau\bv+\bw)],$$
and 
$$\bu(i,j,k)=(\frac{i}{d_1}-\frac{i+j+k}{d_0},\frac{i+j+k}{d_0}-\frac{k}{d_2}).$$
It remains to collect terms with a given basis element $e_{03}(l)$.
\ed

\subsection{Proof of Theorem \ref{main-thm1}}

Let us first fix $\tau$ and assume that $v_1$ and $v_2$ are sufficiently generic
real numbers.
Then Lemmata \ref{connection-lem} and \ref{triple-prod-lem} imply that
for a Massey system $\bc$ one has

\begin{equation}\label{MP-F-eq}
MP(\bc)_l=\exp(\pi i\tau Q(\bv))\sum_{ijk}c_{ijk}\Th_{ijkl}.
\end{equation}

We are going to rewrite this expression in the form similar to $\Th_{Q,f}$. 
As in section \ref{double-sec},
it is convenient to extend the function $(i,j,k)\mapsto c_{ijk}$ by zero
from $\Z^3/\bd\Z^3$ to $\Q^3/\bd\Z^3$ 
(recall that $\bd\Z^3=d_1\Z\times d\Z\times d_1\Z$).
Let us fix $l\in\Z/d_0\Z$ and consider the map
$$\phi_l:\Q^2\to\Q^3/\bd\Z^3:
(\la_1,\la_2)\mapsto (d_1\la_1+\frac{d_1l}{d_0},d_2\la_2-d_1\la_1-\frac{dl}{d_0},
-d_2\la_2+\frac{d_2l}{d_0}).$$
It is easy to check that we have
$\bla\in\La^0(i,j,k,l)$ iff $\phi_l(\bla)=(i,j,k)$.
Hence, we can rewrite \eqref{MP-F-eq} as follows:
\begin{equation}\label{MP-main-eq}
MP(\bc)_l=\exp(\pi i\tau Q(\bv))
\sum_{\bla\in\Q^2\cap (C-\bv)}
\sign(\la_1+v_1)c_{\phi_l(\bla)}
\exp[\pi i\tau Q(\bla)+2\pi i\bla\cdot(\tau\bv+\bw)],
\end{equation}
where $\bla=(\la_1,\la_2)$. We claim that the sums of the function 
$$\bla\mapsto c_{\phi_l(\bla)} 
\exp[\pi i\tau Q(\bla)+2\pi i\bla\cdot(\tau\bv+\bw)]$$
over vertical and horizontal lines in $\Q^2$ are zero. This
is equivalent to the following system of identities:
\begin{equation}\label{vanish-eq}
\sum_{ijk}c_{ijk}\sum_{\bla\in\La^0(i,j,k,l)\cap  L_r(a)}
\exp[\pi i\tau Q(\bla)+2\pi i\bla\cdot(\tau\bv+\bw)]=0,\ a\in\Q,\ r=1,2,
\end{equation}
where $L_r(a)=\{(\la_1,\la_2):\la_r=a\}$.
To prove the identity corresponding to $r=2$ let us change the variables
$(\la_1,\la_2)$ to $(m_1,\la_2)$ where
$$m_1=\frac{d_2\la_2+(d-d_1)\la_1}{d}.$$
It is easy to check that in these new variables the quadratic form $Q$ takes the 
form
$$Q(m_1,\la_2)=\frac{d_1d}{d-d_1}m_1^2-\frac{d_0d}{d-d_1}\la_2^2.$$
On the other hand, the conditions defining $\La^0(i,j,k,l)$ become
\begin{align*}
\La^0(i,j,k,l)=\{(m_1,\la_2):\ 
&m_1\equiv\frac{i}{d_1}+\frac{j}{d}\mod(\Z),\\
&\la_2\equiv\frac{l}{d_0}-\frac{k}{d_2}\mod(\Z),\\
&\frac{dm_1-d_0\la_2}{d-d_1}\equiv\frac{i}{d_1}-\frac{k}{d_2}\mod(\Z).
\end{align*}
Hence, the set $\La^0(i,j,k,l)\cap L_2(a)$ is empty unless 
$$a\equiv\frac{l}{d_0}-\frac{k}{d_2}\mod(\Z).$$
This congruence implies that 
$$d_0a\equiv-\frac{d_0k}{d_2}\equiv\frac{(d-d_1)k}{d_2}\mod(\Z).$$
Hence, there exists $u\in\Z/(d-d_1)\Z$ such
that
$$\frac{u}{d-d_1}\equiv\frac{d_0a}{d-d_1}-\frac{k}{d_2}\mod(\Z).$$
In this situation we have
$$\La^0(i,j,k,l)\cap L_2(a)=\{(m_1,a):\ -m_1\in I_1(u,i,j)\}.$$
Now using the above formula for $Q(m_1,\la_2)$ one can easily verify that 
for fixed $k\in\Z/d_2\Z$ the identity
$$\sum_{ij}c_{ijk}\sum_{\bla\in\La^0(i,j,k,l)\cap  L_2(a)}
\exp[\pi i\tau Q(\bla)+2\pi i\bla\cdot(\tau\bv+\bw)]=0$$
is either trivial (when $\La^0(i,j,k,l)=\emptyset$ for all $(i,j)$) or
equivalent to \eqref{double-eq1} with $u$ chosen as above. 
Summing over all $k$ we derive \eqref{vanish-eq} for $r=2$.
The proof for $r=1$ is very similar: one should start by changing the variables 
$(\la_1,\la_2)$ to $(\la_1,m_2)$, where
$$m_2=\frac{d_1\la_1+(d-d_2)\la_2}{d}.$$
It follows that the summation pattern in the sum
\eqref{MP-main-eq} can be replaced with the pattern used to define indefinite 
theta series (see section \ref{indef-theta-sec}):
\begin{equation}\label{MP-main2-eq}
\exp(-\pi i\tau Q(\bv))MP(\bc)_l=
(\sum_{\la_1\ge 0,\la_2\ge 0}-\sum_{\la_1<0,\la_2<0})
c_{\phi_l(\bla)}\exp[\pi i\tau Q(\bla)+2\pi i\bla\cdot(\tau\bv+\bw)].
\end{equation}
Now we claim that this formula holds also for arbitrary $v_1$ and $v_2$
(not necessarily generic). Indeed, by Lemma \ref{Mas-cont-lem} the
Massey product $MP(\bc)_l$ is a continuous function of $(\bc,v_1,v_2)$
varying in a total space of a vector bundle over $\R^2$. Since the
RHS of \eqref{MP-main2-eq} also depends continuously on $(\bc,v_1,v_2)$
our claim follows.

Finally, let us consider the situation of Theorem \ref{main-thm1},
so we assume that $\bc$ is a universal Massey system
(and $v_1$ and $v_2$ are rational). Making the change of variables 
$\bx=\bla+\bv$ and using \eqref{vanish-eq}
we can rewrite \eqref{MP-main2-eq} as follows:
\begin{align*}
&MP(\bc)_l=\exp(-2\pi i\bv\cdot\bw)
(\sum_{x_1\ge 0,x_2\ge 0}-\sum_{x_1<0,x_2<0})
c_{\phi_l(\bx-\bv)}\exp(\pi i\tau Q(\bx)+2\pi i\bx\cdot\bw)=\\
&(\sum_{x_1\ge 0,x_2\ge 0}-\sum_{x_1<0,x_2<0})f(\bx)q^{Q(\bx)/2},
\end{align*}
where for $\bx=(x_1,x_2)\in\Q^2$ we set
\begin{equation}\label{f-eq}
f(\bx)=f_{\bc,l}(\bx)=c_{\phi_l(\bx-\bv)}\exp(2\pi i(\bx-\bv)\cdot\bw).
\end{equation}
As we have seen above the sums of $f(\bx)q^{Q(\bx)/2}$ over all vertical and
horizontal lines in $\Q^2$ are zero. 
Hence, $f$ satisfies \eqref{f-condition} and we are done.
\ed

\subsection{Proof of Theorem \ref{main-thm2}}

We start by analyzing in more detail
the correspondence $\bc\mapsto f_{\bc,l}$ constructed 
in the proof of Theorem \ref{main-thm1}.
In what follows we assume that the parameters $\bv$ and $\bw$ are equal to zero.

It is convenient to change the meaning of the parameter $l$ by allowing it to run
through $\Q$, so that our previous expressions containing $l$ should be
viewed as $d_0\Z$-periodic functions of $l$ supported on $\Z$.
Let us consider a homomorphism
$$\phi:Q^3\to\Q^3:
(\la_1,\la_2,l)\mapsto (d_1\la_1+\frac{d_1l}{d_0},d_2\la_2-d_1\la_1-
\frac{dl}{d_0},
-d_2\la_2+\frac{d_2l}{d_0}).$$
For $l\in\Z$ it is related to the homomorphisms $\phi_l$ introduced in
the proof of Theorem \ref{main-thm1}:
$$\phi_l(\la_1,\la_2)=\phi(\la_1,\la_2,l)\mod \bd\Z^3.$$
It is easy to see that $\phi$ is invertible and that
$$\phi^{-1}(i,j,k)=
(\frac{i}{d_1}-\frac{i+j+k}{d_0},-\frac{k}{d_2}+\frac{i+j+k}{d_0},i+j+k).$$
Let us set $\wt{\Ga}:=\phi^{-1}(\Z^3)$. The above formula for $\phi^{-1}$ shows 
that 
$\wt{\Ga}$ is contained in $\Q^2\times\Z$, so we can write
$$\wt{\Ga}=\sqcup_{l\in\Z}\Ga(l)\times\{l\},$$
with some $\Ga(l)\sub\Q^2$.
Since the restriction of $\phi$ to $\Q^2\times\{l\}$ is essentially $\phi_l$,
we obtain that
$$\Ga(l)=\cup_{(i,j,k)\in\Z^3}\La^0(i,j,k,l).$$
Let us also consider the lattice 
$$\Ga=p_{12}(\wt{\Ga})=\cup_{l\in\Z}\Ga(l),$$
where $p_{12}:\Q^3\to\Q^2$ is the projection on the first two components.
When we need to show the dependence of $\Ga$ on $(d_1,d_2,d)$
we will write $\Ga_{d_1,d_2,d}$.
From the relation $\Ga=p_{12}\phi^{-1}(\Z^3)$ it is clear 
that for any positive rational number $N$ we have 
$$\Ga_{Nd_1,Nd_2,Nd}=\frac{1}{N}\Ga_{d_1,d_2,d}.$$

Let us set $\wt{\La}:=\phi^{-1}(\bd\Z^3)$.
One can easily check that
$$p(\wt{\La})=\La,$$
where $\La$ was defined in section \ref{triple-sec}.
It will turn out to be crucial for us that unlike $\Ga$,
the lattice $\La$ does not change when we
rescale $(d_1,d_2,d)$ to $(Nd_1,Nd_2,Nd)$.

Set $\be_3=(0,0,1)\in\Q^3$.
Note that we have the inclusions $d_0\Z \be_3\sub\wt{\La}\sub\wt{\Ga}$.
Let us denote by $p:\wt{\Ga}/d_0\Z \be_3\to\Ga$ 
the map induced by the projection $p_{12}$.
Then $p$ is surjective with the kernel 
$$\phi^{-1}(\Z^3)\cap \Q \be_3/d_0\Z \be_3=\frac{d_0}{g}\Z \be_3/d_0\Z \be_3,$$ 
where $g$ is the greatest common divisor of $d$, $d_1$ and $d_2$. 
As before we view $\bc=(c_{ijk})$ as a 
$\bd\Z^3$-periodic function on $\Q^3$ supported on
$\Z^3$. For such a function we have
$$f_{\bc}:=\sum_{l\in\Z/d_0\Z}f_{\bc,l}=p_!\phi^*\bc,$$
where $\phi^*$ denotes the pull-back and $p_!$ the push-forward (the summation
over fibers of $p$). 

The homomorphisms $\phi$ and $p_{12}$ induce isomorphisms
\begin{equation}\label{main-isom-eq}
\Z^3/(\bd\Z^3+\frac{d_0}{g}\Z\phi(\be_3))\simeq 
\wt{\Ga}/(\frac{d_0}{g}\Z\be_3+\wt{\La})\simeq\Ga/\La,
\end{equation}
so we get a natural identification between the space of 
$(\bd\Z^3+\frac{d_0}{g}\Z\phi(\be_3))$-periodic functions on $\Z^3$ and 
the space of $\La$-periodic functions on $\Ga$.
Hence, if $\bc$ is $\frac{d_0}{g}\Z\phi(\be_3)$-periodic 
then $\frac{1}{g}f_{\bc}$ is
the $\La$-periodic function on $\Ga$ corresponding to $\bc$ under the above
identification. Furthermore, as we have seen in the proof of Theorem 
\ref{main-thm1},
if $\bc$ is a universal Massey system then $f_{\bc}$ satisfies 
\eqref{f-condition}.
The following lemma asserts that the converse is also true for $\bc$ satisfying 
some 
additional periodicity conditions.

\begin{lem}\label{isom-lem} Assume that
$$\frac{2d_1g_1}{d-d_1}\in\Z \text{ and } \frac{2d_2g_2}{d-d_2}\in\Z,$$
where $g_1=gcd(d_1,d)$ and $g_2=gcd(d_2,d)$,
and let us define the following elements in $\wt{\Ga}$:
$$\bz_1=(-\frac{2d_2g_1}{d_0(d-d_1)},\frac{2g_1}{d_0},2g_1),
\ \bz_2=(-\frac{2g_2}{d_0},\frac{2d_1g_2}{d_0(d-d_2)},2g_2).$$
We denote by $\De\sub\wt{\Ga}$ the subgroup generated by 
$\frac{d_0}{g}\be_3$, $\bz_1$ and $\bz_2$.
Let $C(\Z^3/(\phi(\De)+\bd\Z^3))$ be the space of $(\phi(\De)+\bd\Z^3)$-periodic 
functions 
$\bc=(c_{ijk})$ on $\Z^3$
such that systems of equations \eqref{double-eq1} and \eqref{double-eq2} 
hold identically in $\tau$.
Let also $F(\Ga/(p_{12}(\De)+\La))$ be the space of $(p_{12}(\De)+\La)$-periodic
functions $f$ on $\Q^2$ with support in $\Ga$ satisfying \eqref{f-condition}. 
Then the map
$$C(\Z^3/\bd\Z^3+\phi(\De))\to F(\Ga/(p_{12}(\De)+\La)):\bc\mapsto f_{\bc}$$
is an isomorphism.
\end{lem}

\Pf . Since the map $\bc\mapsto \frac{1}{g}f_{\bc}$ is induced by the isomorphism
\eqref{main-isom-eq}, it is clear that it transforms $\phi(\De)$-periodicity
to $p_{12}(\De)$-periodicity. Therefore,
we only have to check that if $\bc$ is a $\phi(\De)+\bd\Z^3$-periodic function on
$\Z^3$ such that $f_{\bc}$ satisfies \eqref{f-condition}, then 
\eqref{double-eq1} and \eqref{double-eq2} hold for $\bc$. 
Now we observe that the periodicity of $\bc$ with respect to the subgroups
$\Z\phi(\bz_1)$ and $\Z\phi(\bz_2)$ is exactly the periodicity given by equations
\eqref{per-eq1} and \eqref{per-eq2}.
Hence, by Lemma \ref{tech-lem} it suffices to check that $\bc$ satisfies 
equations \eqref{vanish-oper-eq1} and \eqref{vanish-oper-eq2}.
Consider the operators $\wt{A}=A\times\id$ and $\wt{B}=B\times\id$ on 
$\Q^3=\Q^2\times\Q$,
where 
$$A=\left(\matrix -1 & -\frac{2d_2}{d-d_1} \\ 0 & 1\endmatrix\right) 
\text{ and } 
B=\left(\matrix 1 & 0 \\ -\frac{2d_1}{d-d_2} & -1\endmatrix\right)$$
are the operators appearing in \eqref{f-condition} for our quadratic form $Q$.
Our assertion follows immediately 
from the fact that the operators on $(i,j,k)$ entering in 
equations \eqref{vanish-oper-eq1} and \eqref{vanish-oper-eq2} are
exactly $\phi\wt{A}\phi^{-1}$ and $\phi\wt{B}\phi^{-1}$.
\ed

Now we are ready to prove Theorem \ref{main-thm2}. Let $Q(m,n)=am^2+2bmn+cn^2$,
where $a$, $b$ and $c$ are positive rational numbers such that $D=b^2-ac>0$.
Note that we can change the variables $(m,n)$ to $(xm,yn)$, where
$x$ and $y$ are positive rational numbers, and change the data $(Q,f)$ accordingly
without changing the series $\Th_{Q,f}(q)$. Hence, we are allowed to rescale
$(a,b,c)$ to $(x^2a,xyb,y^2c)$. Since $b^2/ac>1$ using such rescaling
we can achieve that $b>a$ and $b>c$. Now let us define positive rational numbers
$(d_1,d_2,d)$ by setting
$$d_1=\frac{D}{b-c},\ d_2=\frac{D}{b-a},\ d=\frac{d_1d_2}{b}=\frac{D^2}{b(b-a)(b-
c)}.$$
Then $Q=Q_{d_1,d_2,d}$ and the inequalities $d>d_1$, $d>d_2$ and 
$d_0=d_1+d_2-d>0$ hold (the last inequality follows from the formula
$d_0=\frac{dD}{d_1d_2}$). Making the change of variables $m=Nm',n=Nn'$ we find
$$\Th_{Q,f}=\Th_{N^2Q,f'}=\Th_{Q_{N^2d_1,N^2d_2,N^2d},f'},$$
where $f'(m',n')=f(Nm',Nn')$. Assume that $f$ is supported on a lattice 
$T\sub\Q^2$
and is $T^0$-periodic for a sublattice $T^0\sub T$. Then $f'$ is supported
on $\frac{1}{N}T$ and is $\frac{1}{N}T^0$-periodic.
If we choose $N$ sufficiently divisible then all the
numbers 
$$N^2d_1,\ N^2d_2,\ N^2d,\ \frac{2d_1g_1}{d-d_1}\text{ and } \frac{2d_2g_2}{d-
d_2}$$ 
will become integers and the following inclusions will hold: 
$$p_{12}(\De)+\La\sub\frac{1}{N}T^0,\ 
\frac{1}{N}T\sub\frac{1}{N^2}\Ga_{d_1,d_2,d}=\Ga_{N^2d_1,N^2d_2,N^2d}.$$
Note that the lattice $p_{12}(\De)+\La$ does not change when we rescale 
$(d_1,d_2,d)$.
Hence, renaming $(N^2d_1,N^2d_2,N^2d,f')$ to $(d_1,d_2,d,f)$
we reduce ourselves to the situation when $Q=Q_{d_1,d_2,d}$, where
$(d_1,d_2,d)$ are integers satisfying the assumptions of Lemma \ref{isom-lem}, and
the function $f$ is $(p_{12}(\De)+\La)$-periodic and is supported on $\Ga$.
Applying Lemma \ref{isom-lem} we find a universal Massey system $\bc$ such that
$f=f_{\bc}$. Finally, as we have shown in the proof of Theorem \ref{main-thm1}, 
one
has
$$\sum_{l\in\Z/d_0\Z}MP(\bc)_l=\sum_{l\in\Z/d_0\Z}\Th_{Q,f_{\bc,l}}(q)=
\Th_{Q,f_{\bc}}(q).$$
\ed

\end{document}